\documentclass[graybox]{svmult}

\usepackage{mathptmx}       
\usepackage{helvet}         
\usepackage{courier}        
\usepackage{type1cm}        

\usepackage{graphicx}        
\usepackage{multicol}        
\usepackage[bottom]{footmisc}

\usepackage{amsmath}
\usepackage{amsfonts}
\usepackage{amssymb}

\def\CC{{\mathbb C}}

\def\E{{\cal E}}

\def\M{{\cal M}}

\def\RR{{\mathbb R}}

\def\ZZ{{\mathbb Z}}

\def\eps{{\varepsilon}}

\def\é{{\'e}}
\def\à{{\`a}}
\def\è{{\`e}}
\def\ù{{\`u}}
\def\ç{{\c{c}}}

\DeclareMathOperator{\expo}{Exp}
\DeclareMathOperator{\Pf}{Pf}

\def\GR{G^{\diamond}}
\def\G1bar{\bar{G}_1}
\def\GD{\mathcal{G}}

\def\L1bar{\bar{L}_1}

\def\Pdimer{P_{\mathrm{dimer}}}
\def\Zdimer{Z_{\mathrm{dimer}}}
\def\Pising{P_{\mathrm{Ising}}}
\def\Zising{Z_{\mathrm{Ising}}}

\newcommand{\ud}{\mathrm{d}} 

\newenvironment{proofsketch}{
 \noindent{\emph{Proof (sketch)}.}}{
}

\begin{document}

\title*{Statistical mechanics on isoradial graphs}
\author{C\'edric Boutillier and B\'eatrice de Tili\`ere}
\institute{C\'edric Boutillier \at Laboratoire de Probabilit\'es et Mod\`eles Al\'eatoires, Universit\'e Pierre et Marie Curie, 
4 place Jussieu, F-75005 Paris. \email{cedric.boutillier@upmc.fr}
\and B\'eatrice de Tili\`ere \at Laboratoire de Probabilit\'es et Mod\`eles Al\'eatoires, Universit\'e Pierre et Marie Curie, 4 place Jussieu, F-75005 Paris. \email{beatrice.de\_tiliere@upmc.fr}}

\maketitle

\abstract{Isoradial graphs are a natural generalization of regular graphs which give, for many models of statistical
mechanics, the right framework for studying models at criticality.  
In this survey paper, we first explain how isoradial graphs naturally arise in two approaches used by physicists:
transfer matrices and conformal field theory. This leads us to the fact that isoradial graphs provide a natural setting
for discrete complex analysis, to which we dedicate one section. Then, we give an overview of explicit results obtained 
for different models of statistical mechanics defined on such graphs: the critical dimer model when the underlying graph
is bipartite, the $2$-dimensional critical Ising model, random walk and spanning
trees and the $q$-state Potts model.}  


\section{Introduction}

Statistical mechanics aims at describing large scale properties of physics systems based on models
which specify interactions on a microscopic level. In this setting, physics systems are modelled 
by random configurations of graphs embedded in a $d$-dimensional space. Since vertices typically represent atoms,
the goal is to let the mesh size tend to zero and rigorously understand 
the limiting behavior of the system. 

Although real world suggests that we should focus on the case of $3$-dimensional systems, 
we restrict ourselves to the case of dimension $2$. There are two main reasons guiding this
choice: first, for models we consider, only very few rigorous
results exist in dimension $3$, and more importantly, $2$-dimensional systems exhibit
beautiful and rich behaviors which are strongly related to this choice of dimension.

Historically, the most studied graph is certainly $\ZZ^2$, followed by the honeycomb and
triangular lattices. Since solving statistical mechanics models involves dealing with combinatorial
and geometric features of the underlying graph, it has been most convenient to handle the 
simplest and most regular ones. However in this paper, we deal with models defined on a 
more general class of graphs, called \emph{isoradial graphs}. 

The motivation behind this generalization is to find the most appropriate setting, 
which exhibits some essential features required to solve the questions addressed, thus allowing for 
proofs revealing the true nature of the problems, and hopefully a full understanding of the issues 
at stake. Before going any further, let us define isoradial graphs and the corresponding rhombus graph.

\begin{definition}
A graph $G=(V,E)$ is said to be {\em isoradial}~\cite{Kenyon3}, 
if it has an embedding in the plane such that every face is inscribed in a circle of radius~1, and such that 
all circumcenters of the faces are in the closure of the faces.
\end{definition}
From now on, when we speak 
of the graph $G$, we mean the graph together with a particular isoradial embedding in the
plane. Examples of isoradial graphs are given in Figure~\ref{fig:isingcrit} below, in particular 
$\ZZ^2$, the honeycomb and triangular lattices are isoradial.

\begin{figure}[ht]
  \begin{center}
    \includegraphics[width=\linewidth]{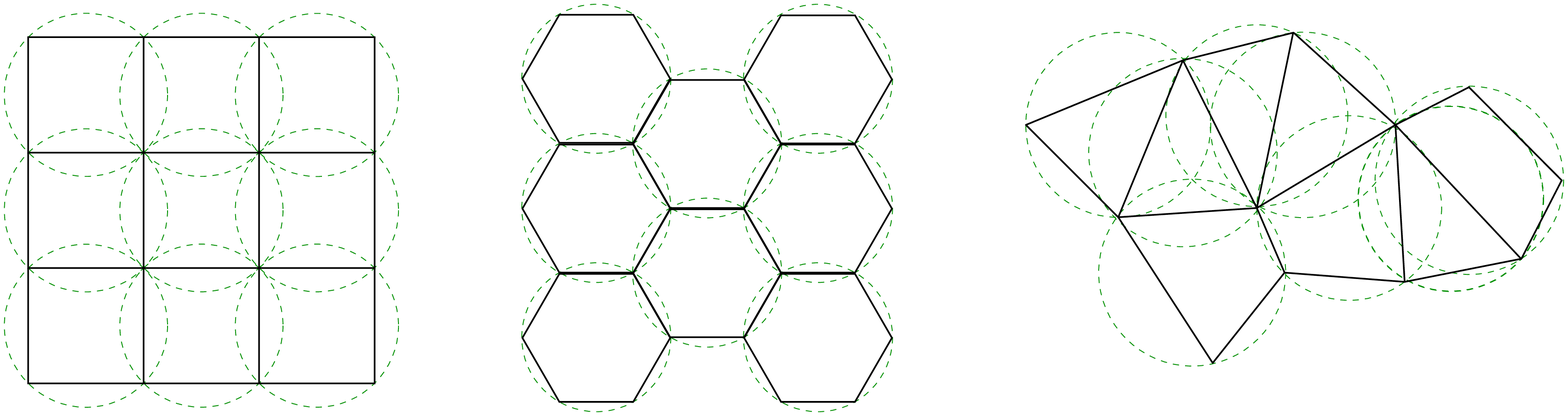}
    \caption{Examples of isoradial graphs: the square lattice (left), the honeycomb lattice (center), and a more generic one (right). Every face is inscribed in a circle of radius 1, represented in dashed lines.}
    \label{fig:isingcrit}
  \end{center}
\end{figure}

\begin{figure}[ht]
  \begin{center}
    \includegraphics[width=\linewidth]{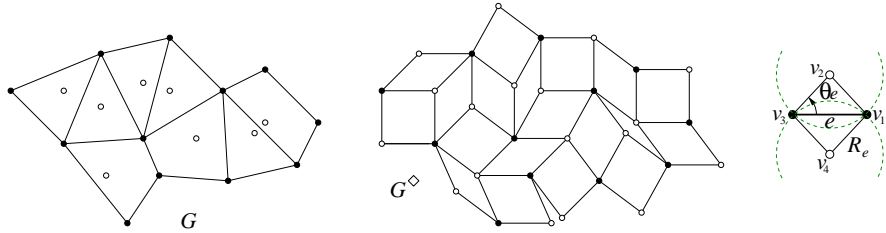}
    \caption{An isoradial graph (left). The white dots are the circumcenters of the faces, which are also the 
vertices of the dual $G^*$. Its diamond graph is represented in the center. On the right is the rhombus $R_e$ and the
half-rhombus angle $\theta_e$ assigned to an edge $e$.}
    \label{fig:iso-rhombi}
  \end{center}
\end{figure}

To such a graph is naturally associated the {\em diamond graph}, denoted by $\GR$, defined as follows. Vertices of $\GR$ consist in the vertices
of $G$, and the circumcenters of the faces of $G$. The circumcenter of each
face is then joined to all vertices which are on the boundary of this face, see
Figure \ref{fig:iso-rhombi}. Since $G$ is isoradial,
all faces of $\GR$ are side-length-$1$ rhombi. Moreover, each edge $e$ of
$G$ is the diagonal of exactly one rhombus $R_e$ of $\GR$; we let $\theta_e$ be the
half-angle of the rhombus at the vertex it has in common with $e$. For later purposes, we label vertices of the 
rhombus $R_e$ by $v_1,\,v_2,\,v_3,\,v_4$, as in Figure~\ref{fig:iso-rhombi} (right).

Isoradial graphs naturally arise in two different approaches to statistical mechanics. We describe them in the
next two sections.

\subsection{Transfer matrices, star transformations and $Z$-invariance}

The first question addressed when solving a model of statistical mechanics is to compute the \emph{free energy},
which is the growth rate of the \emph{partition function}, counting the weighted number of 
configurations. An important technique introduced by physicists to solve two-dimensional models is the use of 
\emph{transfer matrices}, which appears for the first time in the work of Kramers and Wannier 
\cite{KramersWannier1,KramersWannier2}.
If a model of statistical mechanics with local interactions is defined on a torus of size $n\times m$, then the partition function $Z_{m,n}$ can be expressed as 
the trace of the $n$-fold product of a \emph{transfer matrix} $T$, whose rows and columns are indexed by the 
configurations of the model in a strip $m\times 1$, and if $\mathcal{C}$ and $\mathcal{C}'$ are two such configurations, then the matrix element $T_{\mathcal{C},\mathcal{C}'}$ is the contribution to the Boltzmann weight of the interactions in $\mathcal{C}$ and $\mathcal{C}'$ are in consecutive strips.
The free energy of the model $f$ can be expressed in terms of the largest eigenvalue of $T$ denoted by $\lambda_1$:
\begin{displaymath}
f=-\lim_{m,n\rightarrow\infty}\frac{1}{mn}\log(Z_{m,n})=
-\lim_{m\rightarrow\infty}\frac{1}{m}\log(\lambda_1).
\end{displaymath}

If the interaction constants on the graph are not homogeneous, then the transfer matrices used to pass from one strip 
to another are different. If we want to be able to compute the free energy, as in the homogeneous case, in terms of 
spectral characteristics of the transfer matrices, then these should commute.

\begin{figure}[ht]
  \begin{center}
    \includegraphics[width=3.6cm]{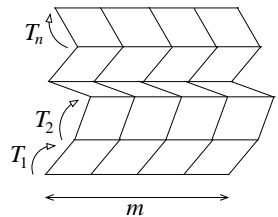}
    \caption{Transfer matrices for an inhomogeneous model on the $m\times n$ torus.}
    \label{fig-transfer1}
  \end{center}
\end{figure}

\begin{figure}[ht]
  \begin{center}
    \includegraphics[width=4cm]{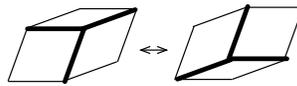}
    \caption{Star-star transformation (the name \emph{star} referring to the $3$-branches stars).}
    \label{fig-transfer2}
  \end{center}
\end{figure}

Since the interactions are local, the transfer matrices are themselves sums of local operators $R$ acting on 
configurations at neighboring sites, and diagrammatically represented by one rhombus of a strip.
A sufficient condition to ensure the commutation of transfer matrices is to demand that these $R$-matrices satisfy 
some algebraic relation, called the \emph{Yang-Baxter equations}, which can be loosely formulated as follows. 
Consider the two ways to tile a hexagon with the same three rhombi on Figure~\ref{fig-transfer2}. The transformation 
from one to the other is called a \emph{star-star transformation}. The model satisfies the Yang-Baxter equations 
if the sum of Boltzmann weights over local configurations of the left hand-side equals those of the right 
hand side. This very strong constraint gives a set of equations which needs to be satisfied by Boltzmann weights,
see \cite{Perk:YangBaxter,Baxter:exactly} for more details.
If a non trivial solution is found, we say that the model is \emph{$Z$-invariant}.

The condition of $Z$-invariance yields a natural generalization. Suppose that the
graph embedded on the torus of size $m\times n$ is a periodic tiling by rhombi (where the period is independent of 
$m$ and $n$). Then
by performing a sequence of star-star moves, this graph can be transformed into big pieces of
tilted copies of $\ZZ^2$ \cite{Kenyon:trieste}. If boundary effects
can be neglected, then it suffices to solve the model on each of the copy of pieces of $\ZZ^2$
using transfer matrices. Summarizing, if the model is $Z$-invariant, then it can be solved on any 
periodic graph consisting of rhombi.

Consider an isoradial graph $G$ and recall that a unique rhombus of the underlying diamond graph $\GR$
can be assigned to every edge of $G$. Then, looking at edges rather than rhombi,
the star-star transformation becomes a \emph{star-triangle 
transformation}, see Figure \ref{fig-transfer4} below. 

\begin{figure}[ht]
  \begin{center}
    \includegraphics[width=4cm]{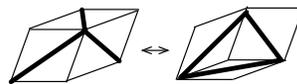}
    \caption{Star-triangle transformation for isoradial graphs.}
    \label{fig-transfer4}
  \end{center}
\end{figure}

The definition of $Z$-invariance naturally extends to this setting, 
and if the model is $Z$-invariant the transfer
matrix approach can be performed on $G^{\diamond}$, thus explaining the occurrence of isoradial graphs and
$Z$-invariance in the context of transfer matrices.

\subsection{Conformal field theory and discrete complex analysis}

Conformal field theory (CFT), introduced by Belavin, Polyakov and Zamolodchikov \cite{BelavinPolyakovZamolodchikov} (see also \cite{diFrancescoMathieuSenechal}) is a theory 
which aims at describing models at criticality, supposed not to depend on specificities of the graph, and to be  
\emph{conformally invariant}.
This very strong statement remained largely inaccessible to the mathematics community (except for a few models, such as dimers on the square lattice \cite{Kenyon:confinv}) until the introduction of the SLE
process by Schramm in 1999. The SLE is conformally invariant and conjectured to be the limiting process for many models
of statistical mechanics, thus filling a huge gap with CFT.
For relations between CFT and SLE, see for example \cite{FriedrichWerner,BauerBernard}. 
Several of these conjectures are now solved: loop erased random walk and
uniform spanning tree (Lawler, Schramm, Werner \cite{LawlerSchrammWerner}), site percolation on the
triangular lattice (Smirnov \cite{Smirnov:perco}), Ising (Smirnov \cite{Smirnov:isingZ2}, Chelkak and Smirnov\cite{ChelkakSmirnov:ising})\dots 

It is thus of key importance to have a setting suitable for doing
discrete complex analysis, and proving convergence to its continuous counterpart. To this purpose, 
isoradial graphs are a perfectly suited object. Indeed, when considering an edge and its dual, they consist in 
the two diagonals of a rhombus, and are thus orthogonal. This allows for a natural discretization of the 
Cauchy-Riemman equation. Since a lot of developments have happened in this direction, we dedicate it the next section.

\section{Discrete complex analysis on isoradial graphs}

As mentioned above, isoradial graphs provide a natural
framework for a generalization of the
construction in the case of the square lattice. The ideas presented here go
back to Duffin \cite{Duffin}, and have been developed later by Mercat \cite{Mercat:ising}, 
Kenyon \cite{Kenyon3}, Chelkak and Smirnov \cite{ChelkakSmirnov:toolbox}, Cimasoni \cite{Cimasoni:dirac} and others.
 There are other possible discretizations of complex analysis on isoradial graphs with remarkable properties. 
 An example using the notion of cross-ratios is given in \cite{BobenkoMercatSuris}.

\subsection{Discrete holomorphic and discrete harmonic functions}\label{sec:21}

Let us start with the notion of discrete holomorphy.  
In the continuous setting, the Cauchy-Riemann equations satisfied by a
holomorphic function imply that its partial derivatives along two orthogonal
unit vectors differ by a factor $i$. In the setting of isoradial graphs, one
takes advantage of the fact that faces of $\GR$ are rhombi, and thus
have perpendicular diagonals, to write down finite difference equations which
are discretizations of Cauchy-Riemann equations. 

Let $f$ be a $\CC$-valued function on vertices of $\GR$, let $e$ be an edge of $G$, and
$R_e$ be the corresponding rhombus with the labeling of its vertices $v_1,\,v_2,\,v_3,\,v_4$, given in 
Figure~\ref{fig:iso-rhombi}.
We say that the function $f$ is \emph{discrete holomorphic} at the rhombus $R_e$, if:
\begin{equation*}
  \frac{f(v_3)-f(v_1)}{ v_3- v_1} = \frac{f(v_4)-f(v_2)}{ v_4 -  v_2}.
\end{equation*}
More generally we define a discrete operator, $\bar\partial:\CC^{V(\GR)}\rightarrow \CC^{E(G)}$, by: 
\begin{equation*}
\bar{\partial} f (e) = i(v_4-v_2)(v_3-v_1)
\left[ \frac{f(v_4)-f(v_2)}{v_4 - v_2} - \frac{f(v_3)-f(v_1)}{v_3 - v_1}\right]
=\sum_{j=1}^4 \bar{\partial}_{e,v_j} f(v_j),
\end{equation*}
where 
\begin{equation}\label{eq:dbar}
\bar{\partial}_{e,v_j}=i(v_{j-1}-v_{j+1}),
\end{equation}
with indices written mod $4$. A function $f$ is then said to be
\emph{discrete holomorphic on $\GR$}, if $\bar\partial f \equiv 0$. Similarly, we define an operator $\partial$
by replacing the coefficients $\bar{\partial}_{e,v_j}$ by their complex conjugate.

The operators $\partial$ and $\bar\partial$ can be extended to 1-forms (from $\CC^{E(G)}$ to $\CC^{V(\GR)}$):
let $g$ be a $\CC$-valued function on edges of $G$, then for all $v\in V(\GR),$
\begin{equation*}
  \quad\partial g (v) = \sum_{e=(v,v')\sim v}\frac{1}{(v'-v)} g(e).
\end{equation*}
The extension of the operator $\bar{\partial}$ is again obtained by replacing coefficients by their complex conjugate.

Note that the normalization of the operators $\partial$ and $\bar{\partial}$ we adopt here differs 
from that chosen by Mercat, or by Chelkak and Smirnov. One can recover their normalization by multiplying 
our operators by diagonal operators acting on edges. So these variations do correspond to the same notion 
of discrete holomorphy.

An notion intimately connected to holomorphy is that of harmonicity. The Laplace operator $\Delta$ on $G$ can be written 
as the restriction of $\partial \bar\partial$ to functions supported on $G$
\footnote{One can define in a similar way a Laplacian 
on $G^*$ by restricting the same operator $\partial\bar{\partial}$ to $G^*$}. Its action on a function $f$ defined 
on vertices of $G$ is given by
\begin{equation*}
  \forall v\in G, \quad \Delta f (v)=\partial(\bar{\partial}f)(v) = \sum_{u\sim v} \tan \theta_{uv} (f(u)-f(v)).
\end{equation*}
A function $f$ is said to be \emph{discrete harmonic} if $\Delta f\equiv 0$. This 
corresponds to choosing for every edge $e$, a conductance equal to $\tan\theta_e$.

From the factorisation $\Delta = \partial \bar\partial$, one sees that, like in the
continuous case, the restriction of a discrete holomorphic function to $G$ (resp. $G^*)$ is 
harmonic. Conversely, if $H$ is a harmonic function on $G$, there exists a
harmonic function $H^*$ on $G^*$ such that $H+i H^*$ is discrete holomorphic.

From this choice of Laplacian, one can define the discrete analogues of classical quantities in continuous potential theory 
(Green function, Poisson kernel, harmonic measure,\dots). A surprising amount of estimates have a discrete version, 
and can be found in \cite{ChelkakSmirnov:toolbox}.

Discrete complex analysis can serve as an approximation of the continuous theory. Mercat \cite{Mercat:ising} proved 
that the pointwise limit of a sequence of discrete holomorphic functions on isoradial graphs in a domain with a mesh 
going to zero is holomorphic. See \cite{Cimasoni:dirac} for a generalisation of this result to isoradial discretization of 
compact manifolds. These local convergence results have been supplemented by global convergence theorems 
\cite{ChelkakSmirnov:ising}, implying in particular the uniform convergence for the discrete potential theory 
objects to their continuous counterparts.

\subsection{Discrete exponential functions}

There is a class of functions playing a special r\^ole in this theory. These functions, called 
\emph{discrete exponential functions}, are defined recursively on $\GR$. For any given vertex $v_0$ of $\GR$ 
and any $\lambda\in \mathbb{C}$, the function $\expo_{v_0}(\cdot; \lambda)$  is defined as follows: 
its value at $v_0$ is 1, and if $v_1$ and $v_2$ are neighbors in $\GR$, let $e^{i\alpha}=v_2-v_1$, then:
\begin{equation*}
  \expo_{v_0}(v_2; \lambda) = \expo_{v_0}(v_1;\lambda) \frac{1 +\lambda e^{i\alpha}}{1-\lambda e^{i\alpha}}.
\end{equation*}


The name comes from the fact that these functions satisfy the following identity:
for any pair of neighbouring vertices $v_1$ and $v_2$ of $\GR$,
\begin{equation*}
  \expo_{v_0}(v_2;\lambda) - \expo_{v_0}(v_1; \lambda) = 
\lambda \frac{\expo_{v_0}(v_2;\lambda) + \expo_{v_0}(v_1; \lambda)}{2}(v_2-v_1),
\end{equation*}
which is a discrete version of the differential equation
\begin{equation*}
\mathrm{d} \exp(\lambda x) = \lambda \exp(\lambda x) \mathrm{d} x.
\end{equation*}
satisfied by the usual exponential function.
It is straightforward to check from the definition that $\expo_{v_0}(\cdot; \lambda)$ is discrete holomorphic. 
In fact, Bobenko, Mercat and Suris \cite{BobenkoMercatSuris} show that any discrete holomorphic function can be written as a (generalized) linear combination of discrete exponential functions, at least  
under the quasicrystallic assumption (a finite number of slopes $\pm e^{i\theta_1},\dots,\pm e^{i\theta_d}$ 
for the rhombi chains):
if $f$ is a discrete holomorphic function on $\GR$, then there exists a function $g(\lambda)$ defined on a neighborhood in $\mathbb{C}$ of $\{\pm e^{i\theta_1},\dots,\pm e^{i\theta_d}\}$  such that
  \begin{equation*}
    \forall v\in \GR, \quad f(v) = \oint_{\Gamma} \expo_{v_0}(v; \lambda) g(\lambda) \mathrm{d}\lambda,
  \end{equation*}
  where $\Gamma$ is a collection of disjoint small positively oriented contours around the possible p\^oles $\pm e^{i\theta_j}$ of the integrand.

\subsection{Geometric integrability of discrete Cauchy-Riemann equations}

An important feature of the discrete Cauchy-Riemann is the property of \emph{3D consistency} \cite{BobenkoMercatSuris}. Let $f$ be a discrete holomorphic function 
on $\GR$. Given its value on three vertices of a rhombus, the other can be computed using the discrete Cauchy-Riemann 
equation. 

If a \emph{star-star} transformation is realised on the graph $\GR$ (see Figure \ref{fig-transfer2}), can we determine the value of the function $f$ at the new vertex in the center, so that $f$ is again discrete holomorphic?
\emph{A priori}, there are three different ways to compute a possible value: one 
for each rhombus. It turns out that the three values obtained are equal.

Under the quasicrystallic assumption, $\GR$ can be seen as a monotone surface in $\mathbb{Z}^d$, projected back to a properly chose plane. The star-star move is a local displacement of this surface along the faces of a cube. If given two monotonic surfaces $\Sigma_1$ and $\Sigma_2$ in $\mathbb{Z}^d$ that can be deduced one from the other by a sequence of star-star moves, then there is a canonical way to extend a discrete holomorphic function $f$ defined on $\Sigma_1$ to $\Sigma_2$ by ``pushing'' its values along the deformation making use of the Cauchy-Riemann equations.
In this sense, this 3D consistency can be considered as a notion of \emph{integrability} \cite{BobenkoSuris}. 
This property is very much related to the \emph{$Z$-invariance} in statistical mechanics.

\subsection{Generalization of the operator $\bar\partial$}\label{sec:24}

The double graph $G^D$ of $G$ is a planar bipartite graph constructed as
follows, see also Figure \ref{fig:double} below. Vertices of $G^D$ are decomposed into two classes: black vertices, which are
vertices of $G$ and of $G^*$, and white vertices, which are the centers of the
rhombi of $\GR$, seen either as edges of $G$ or as those of $G^*$. A black
vertex $b$ and a white vertex $w$ are connected by an edge in $G^D$ if the
vertex corresponding to $b$ is incident to the edge corresponding to $w$ in
either $G$ or $G^*$. That is, edges of $G^D$ correspond to half diagonals of rhombi of $\GR$.
It turns out that $G^D$ has an isoradial embedding (with
rhombi of side-length $\frac{1}{2}$) obtained by splitting each rhombus of $\GR$
into four identical smaller rhombi.

\begin{figure}[ht]
  \begin{center}
    \includegraphics[width=5cm]{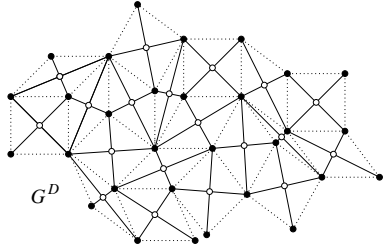}
    \caption{The double graph $G^D$ of the graph on the right of Figure \ref{fig:isingcrit}: in thick full lines is the 
original graph $G$, in full lines is the dual graph $G^*$, in dotted lines is the underlying $\GR$ graph.}
    \label{fig:double}
  \end{center}
\end{figure}

The operator $\bar\partial:\CC^{V(\GR)}\rightarrow \CC^{E(G)}$ of Section \ref{sec:21}
can be interpreted now as an operator $\bar\partial:\CC^{B(G^D)}\rightarrow \CC^{W(G^D)}$, where $B(G^D)$ (resp.
$W(G^D)$) denotes the set of black (resp. white) vertices of $G^D$. Let $f$ be a $\CC$-valued function defined on 
black vertices of $G^D$, then for every white vertex of $G^D$:
\begin{equation*}
\bar{\partial} f (w)=\sum_{b\sim w}\bar{\partial}_{w,b}f(b).
\end{equation*}
Let $w,x,b,y$ be the vertices of a rhombus enumerated in cclw order, so that $w\in W(G^D)$ and $b\in B(G^D)$. Then by Formula \eqref{eq:dbar}, 
\begin{equation}\label{eq:dbar2}
\bar{\partial}_{w,b}=i(x-y).
\end{equation}
This definition can be extended to any bipartite isoradial graph, see Section \ref{sec:32}.

For a detailed study
of this operator on isoradial graphs on compact surfaces and its connections
with discretization of geometric structures (especially spin structures), see
\cite{Cimasoni:dirac}.

\section{Dimer model}

The \emph{dimer model} represents the adsorption of 
diatomic molecules on the surface of a crystal. The surface of the crystal is modeled by
a planar graph $G=(V,E)$, which we assume to be finite for the moment. 
A \emph{dimer configuration} of $G$ is a \emph{perfect matching}
of $G$, that is a subset of edges $M$ such that every vertex of $G$ is incident to a unique 
edge of $M$. Denote by $\M(G)$ the set of dimer configurations of $G$. An example of
dimer configuration is given in Figure \ref{fig:Fisher_correspondence}. Assume that a 
positive \emph{weight function} $\nu$ is assigned to edges of $G$, that is every edge $e$ of $G$ has
weight $\nu_e$. Then, when the graph $G$ is finite, the probability of occurrence of a dimer 
configuration $M$ is given by the \emph{dimer Boltzmann measure}:
\begin{equation*}
\Pdimer(M)=\frac{\prod_{e\in M}\nu_e}{\Zdimer},
\end{equation*}
where $\Zdimer=\sum_{M\in\M(G)}\prod_{e\in M}\nu_e$ is the normalizing constant known as the 
\emph{partition function}.

The dimer model has the attractive feature of being exactly solvable, meaning that there is
an explicit expression for the partition function. This fundamental result is due to 
Kasteleyn \cite{Kast61} and independently to Temperley and Fisher \cite{TemperleyFisher}. It relies on the 
Kasteleyn matrix, which we now define. A 
\emph{Kasteleyn orientation of the graph $G$} is an orientation of the edges such that every 
face is \emph{clockwise odd} meaning that, when traveling clockwise around edges of
a face, an odd number of them are co-oriented. The corresponding \emph{Kasteleyn matrix} $K$ is the 
associated signed, weighted, adjacency matrix of the graph: lines and columns of $K$ are 
indexed by vertices of $G$, and for any two vertices $x,y$ of $G$, the coefficient 
$K_{x,y}$ is:
$$
K_{x,y}=
\begin{cases}
\nu_{xy}&\text{ if }x\sim y\text{ and the edge $xy$ is oriented from $x$ to $y$}\\
-\nu_{xy} &\text{ if }x\sim y\text{ and the edge $xy$ is oriented from $y$ to $x$}\\ 
0&\text{ otherwise}.
\end{cases}
$$
\begin{theorem}[\cite{Kast61,TemperleyFisher}]\label{thm:KTF}
The dimer partition function is, 
$$
\Zdimer=|\Pf(K)|=\sqrt{\det(K)}.
$$
\end{theorem}
\begin{remark} 
The Pfaffian, denoted by $\Pf$, of a skew-symmetric matrix is a polynomial in the entries of the matrix, which is a square root of the determinant.
\end{remark}
\begin{proofsketch}
Refer to \cite{Kasteleyn} for details. When writing out the Pfaffian of an adjacency matrix as a sum over permutations,
there is exactly one term per dimer configuration. The issue is that each term
comes with a signature. The purpose of the Kasteleyn orientation is to
compensate the signature with signs of coefficients.
\end{proofsketch}

As a consequence of Theorem \ref{thm:KTF}, Kenyon derives an explicit
expression for the dimer Boltzmann measure.
\begin{theorem}[\cite{Ke:LocStat}]
Let $E=\{e_1=x_1 y_1,\cdots,e_k=x_k y_k\}$ be a subset of edges of the graph $G$. Then, the probability
that these edges are in a dimer configuration chosen with respect to the dimer Boltzmann measure is;
\begin{equation*}
\Pdimer(\{e_1,\cdots,e_k\})=\left|
\left(\prod_{i=1}^k K_{x_i,y_i}\right)
\Pf((K^{-1})_E)\right|, 
\end{equation*}
where $(K^{-1})_E$ is the submatrix of $K^{-1}$ whose lines and columns are indexed by vertices of $E$.
\end{theorem}
\begin{proof}
Since the proof is very short, we repeat it here. The probability of dimer configurations
containing edges of $E$ is given by the weighted number of these configurations divided by the weighted
number of all configurations. By expanding the Pfaffian along lines and columns, this yields:
\begin{equation*}
\Pdimer(\{e_1,\cdots,e_k\})=\left|\left(\prod_{i=1}^k K_{x_i,y_i}\right)\frac{\Pf(K_{E^c})}{\Pf(K)}\right|.
\end{equation*}
Using Jacobi's formula for Pfaffians \cite{ishikawa}, $|\Pf(K_{E^c})|=|\Pf(K)\Pf((K^{-1})_{E})|$ yields the
result.
\end{proof}

In Section \ref{sec:32}, we state Kenyon's results for the $\bar{\partial}$ operator on infinite, bipartite, isoradial
graphs. Then, in Section \ref{sec:33}, we relate them to the corresponding dimer model.

\subsection{Dirac operator and its inverse}\label{sec:32}

As mentioned in Section \ref{sec:24}, the $\bar{\partial}$ operator, also referred to as the \emph{Dirac operator} in
\cite{Kenyon3}, can be extended to infinite bipartite isoradial graphs. Let $G$ be such a graph,
the set of vertices can be divided into two subsets $B\cup W$, 
where vertices in $B$ (the black ones) are only incident to vertices in $W$ (the white ones). Then
$\bar{\partial}$ maps $\CC^B$ to $\CC^W$: let $f$ be a $\CC$-valued function on black vertices, then
$$
(\bar{\partial}f)(w)=\sum_{b\sim w}\bar{\partial}_{w,b}f(b),
$$
where, if $w,x,b,y$, are the vertices of the rhombus $R_{wb}$ in cclw order,
$\bar{\partial}_{w,b}=i(x-y)$ is given by Formula \eqref{eq:dbar2}.

One of the main results of the paper \cite{Kenyon3} is an explicit expression for the coefficients of 
the inverse $\bar{\partial}$, as a contour integral of an integrand 
which only depends on a path joining the two vertices. In order to state the result, let us introduce the following
functions defined on vertices of $G$. Fix a reference white vertex $w_0$ and let $\lambda$ be a complex parameter. Let 
$v$ be a vertex of $G$, and $w_0=v_0,v_1,\cdots,v_k=v$ be a path in $\GR$ from $w_0$ to $v$. Each edge $v_jv_{j+1}$
has exactly one edge of $G$ (the other is an edge of $G^*$). Direct it away from this vertex if it is white, and towards it
if it is black, and let $e^{i\alpha_j}$ be the corresponding vector. The function $f_v$ is defined inductively along the path,
starting from $f_{v_0}\equiv 1$, and 
$$
f_{v_{j+1}}(\lambda)=
\begin{cases}
f_{v_j}(\lambda)(\lambda-e^{i\alpha_j})&\text{ if the edge $v_j v_{j+1}$ leads away from a black or towards}\\
&\text{ a white vertex},\\
\frac{f_{v_j}(\lambda)}{(\lambda-e^{i\alpha_j})}&\text{ otherwise}.
\end{cases}
$$
It is easy to see that the function is well defined, i.e. independent of
the choice of path from $w_0$ to $v$. An important point is that, by using a branched cover of the plane over $w_0$, 
the angles $\alpha_j$ are defined in $\RR$ and not $[0,2\pi[$. Then, Kenyon has the following theorem.
\begin{theorem}\label{thm:dirac}
Coefficients of the inverse Dirac operator can explicitly be expressed as:
\begin{equation}\label{eq:thmdirac}
(\bar{\partial}^{-1})_{b,w_0}=\frac{1}{4\pi^2 i}\oint_{C}f_b(\lambda)\log(\lambda)\ud \lambda,
\end{equation}
where $C$ is a closed contour oriented counterclockwise, containing all poles of the integrand and with the
origin in its exterior.
\end{theorem}

\begin{remark}
The remarkable features of this theorem are the following.
Explicit computations become tractable, since they only involve residues of rational functions. Moreover, the Formula 
\eqref{eq:thmdirac} has the surprising feature of being ``local'', meaning that if the graph is changed away 
from the vertices $w_0$, $b$ and a path joining them, the value of the inverse Dirac operator stays the same.
\end{remark} 

\subsection{Dirac operator and dimer model}\label{sec:33}

The Dirac operator $\bar{\partial}$ can be represented by an infinite matrix $K$, whose lines (resp. columns) 
are indexed by white (resp. black) vertices of $G$, and the coefficient $K_{w,b}$ is given by:
\begin{equation*}
K_{w,b}=\bar{\partial}_{w,b}=i(x-y)=2\sin(\theta_{wb})e^{i\alpha_{wb}}, 
\end{equation*}
where $\theta_{wb}$ is the half-angle of the rhombus $R_{wb}$, and $e^{i\alpha_{wb}}$ is the unit-vector
in the direction from $w$ to $b$. The matrix $K$ resembles a Kasteleyn matrix of the graph $G$, but differs 
in two instances.
\begin{enumerate}
\item Rows and columns are indexed by only ``half'' the vertices. To overcome this, one can define a weighted
  adjacency matrix $\tilde{K}$ whose lines and columns
are indexed by all vertices of the graph: 
$\tilde{K}=\bigl(\begin{smallmatrix}
0&K\\
-{K}^T&0        
\end{smallmatrix}
\bigr)$. Then, we would have
$\Pf(\tilde{K})=\pm\sqrt{\det(\tilde{K})}=\pm|\det(K)|$. 
\item The weights of the edges have an extra complex factor of modulus one, instead of
a sign given by a Kasteleyn orientation. In \cite{Kenyon3}, Kenyon shows that the matrix is \emph{Kasteleyn flat}, meaning that for
each face of $G$ whose vertices are labeled by $u_1,v_1,\cdots,u_m,v_m$, the coefficients of the matrix $K$ satisfy:
\begin{equation*}
\arg[K_{u_1,v_1}\cdots K_{u_m,v_m}]=\arg[(-1)^{m-1}K_{v_1,u_2}\cdots K_{v_{m},u_{1}}].
\end{equation*}
Using a result of Kuperberg \cite{Kuperberg}, this implies that when the graph is finite and planar, $K$ behaves as a usual Kasteleyn
matrix, thus yielding an explicit formula for the partition function.
\end{enumerate}
One of the main conjectures of the paper \cite{Kenyon3} is to show that the inverse $\bar{\partial}$ operator, 
which we write $K^{-1}$ using the matrix notation, is related to the dimer model on the infinite, bipartite, isoradial, 
graph $G$. Its relation to the Gibbs measure is given in the following.

\begin{theorem}[\cite{Bea1}]\label{thm:gibbsdimer}
There is a unique Gibbs measure $\Pdimer$ on dimer configurations of $G$, such that the probability
of occurrence of a subset of edges $E\{e_1=w_1 b_1,\cdots,e_k=w_k b_k\}$ in a dimer configuration of $G$ chosen
with respect to the Gibbs measure $\Pdimer$ is:
\begin{equation}\label{eq:gibbs}
\Pdimer(\{e_1,\cdots,e_k\})=\left(\prod_{i=1}^k K_{w_i,b_i}\right)\det_{1\leq i,j\leq k} [(K^{-1})_{b_i,w_j}],
\end{equation}
where $(K^{-1})_{b_i,w_j}=(\bar{\partial}^{-1})_{b_i,w_j}$ is given by \eqref{eq:thmdirac}.
\end{theorem}
\begin{proofsketch}
Fix the edge set $E$, and cut out a finite piece of the rhombus graph $\GR$ containing $E$ and paths joining vertices
of $E$. Use a result of \cite{Bea1} by which any finite piece of the graph $\GR$ can be filled with
rhombi in order to become part of a periodic rhombus tiling of the plane. Define the probability of the edge
set $E$ as the weak limit of the Boltzmann measures on the natural toroidal exhaustion of the periodic graph. 
Use the uniqueness of the inverse Dirac operator and its locality property to deduce that this expression coincides with \eqref{eq:gibbs}.
Use Kolmogorov's extension theorem to show existence of a unique measure on $G$ which specifies as \eqref{eq:gibbs} 
on cylinder events.
\end{proofsketch}

\subsection{Other results}

\textbf{Free energy}. Assume that the infinite bipartite isoradial graph is also periodic, and let $G_1=G/\ZZ^2$.
In \cite{Kenyon3}, Kenyon proves an explicit expression for the ``logarithm of the normalized log of the Dirac operator'', which 
has the surprising feature of only depending on rhombus angles of $G_1$:
$$
\log({\det}_1\bar{\partial})=\frac{1}{|V(G_1)|}\sum_{e\in E(G_1)}\left(\frac{1}{\pi}L(\theta)+\frac{\theta}{\pi}
\log(2\sin\theta_e)\right).
$$
Kenyon conjectures it to be related to the \emph{free energy} of the corresponding dimer model, result that we prove this conjecture
in \cite{Bea}.\\

\noindent\textbf{Gaussian free field}.
There is a natural way of assigning a height function to dimer configurations of $G$. In \cite{BeaGFF}, we show that when dimer 
configurations are chosen with respect to the Gibbs measure of Theorem \ref{thm:gibbsdimer}, the fluctuations 
of the height function are described by a Gaussian Free Field. Note that this proof holds for all dimer models
defined on bipartite isoradial graphs, in particular for $\ZZ^2$ and the honeycomb lattice.

\section{Ising Model}

The Ising model, first considered by Lenz \cite{Lenz}, has been introduced in
the physics literature as a model for ferromagnetism. The vertices of a graph
$G=(V,E)$ represent atoms in a crystal with magnetic moment reduced to one component
that can take two values $\pm 1$. A \emph{spin configuration} is thus a function
$\sigma$ on $V$, with values in $\{+1,-1\}$.  Magnets with opposite moments tend
to repel each other, which has a cost in terms of energy.  From these
considerations, the \emph{energy} of a configuration $\sigma$ on a finite graph is defined by:
\begin{equation*}
  \mathcal{E}(\sigma) = - \sum_{e=uv\in E} J_e \sigma_u \sigma_v, 
\end{equation*}
where the positive numbers $(J_e)$ are called \emph{interaction constants}. The
probability of a occurrence of a spin-configuration is then defined using the \emph{Ising Boltzmann measure}:
\begin{equation*}
  \Pising(\sigma)=\frac{\exp(-\mathcal{E}(\sigma))}{\Zising}.
\end{equation*}
The normalising factor, $\Zising=\sum_{\sigma}\exp(-\mathcal{E}(\sigma))$ is the
\emph{Ising partition function} on~$G$.

If the graph is drawn on an orientable surface with no boundary\footnote{The
correspondence can be extended to surfaces with boundary by including in
addition to the closed contours a certain number of paths connected to
boundary.}, then the partition function can be written as a combinatorial sum
over contours. There are in fact two well-known such expansions: the so-called
\emph{high-temperature expansion} (contours on $G$) and the
\emph{low-temperature expansion} (contours on $G^*$). The contours in the
low-temperature expansion have a nice interpretation in terms of the dual
closed curves separating zones of different spins.

The most studied case is the Ising model on (pieces of) regular lattices where
the interaction constants are taken to be constant equal to $\beta$, the
inverse temperature. Kramers and Wannier \cite{KramersWannier1} discovered a duality
for the Ising model on the square lattice:
the measure on contours on $G^*$ coming from the low-temperature expansion is
equal to that obtained by considering the high temperature expansion
of the Ising model on $G^*$ at another temperature $\beta^*$ satisfying $\sinh(2\beta)\sinh(2\beta^*)=1$. This showed that the self-dual temperature for which $\beta=\beta^*=\log\sqrt{1+\sqrt{2}}$ is the critical point if this one is unique. Later, Lebowitz and Martin-L\"of \cite{LebowitzMartinLoef} proved that this is indeed the case.

Note that the case of $\mathbb{Z}^2$ is singular since the graph is isomorphic
to its dual. This duality can be generalized to any graph and its dual. For
example the low-temperature contour measure for the Ising model on the
honeycomb lattice with inverse temperature $\beta$ is the same as the
high-temperature measure for the Ising model on its dual which is the
triangular lattice. But using star-triangle transformations, one can go back to
the Ising model on the honeycomb lattice. This generalizes the duality we had for
the square lattice.

On isoradial graphs, instead of taking the same interaction constants for all
edges, it is natural to take $J_e$ to be a function of the angle $\theta_e$.
In order to ensure some integrability at the discrete level, we can impose that
the Boltzmann weights satisfy the star-triangle relations, yielding a
one-parameter family of $Z$-invariant interaction constants.
There is a generalized duality relation: the
high temperature expansion on $G$ corresponds to the low temperature expansion
on $G^*$ for a dual value of parameter.  We qualify as \emph{critical} the interaction
constants corresponding to the self-dual value of the parameter, given by the following formula,

\begin{equation*} 
  J(\theta_e)=\frac{1}{2}\log\left(
  \frac{1+\sin\theta_e}{\cos\theta_e} \right),
\end{equation*}
and refer to the corresponding Ising model as \emph{critical $Z$-invariant}.
See \cite{Baxter:exactly}, chapter 7.13 for the parametrization of the star-triangle relation.
They coincide with the critical inverse temperatures of the homogeneous triangular, square, and honeycomb 
lattices for $\theta_e$ identically equal to $\frac{\pi}{6}$, $\frac{\pi}{4}$, $\frac{\pi}{3}$ respectively.
The Ising model on two-dimensional graphs is in correspondence with other
well-known models of statistical mechanics: the dimer model and the $q$-random
cluster model with $q=2$.  We differ the discussion of the correspondence with
dimers to Subsection \ref{subseq:Ising-dimers}.

There are essentially two ways of taking the large graph limit of the Ising
model. Either you let the mesh tend to zero at the same time as the number of
vertices of the graph goes to infinity, in order to get continuous objects in a
bounded region of the plane (macroscopic level); or you let the graph tend to
infinity, keeping the mesh size fixed, to get a model of statistical mechanics
on an infinite graph.

The first approach has been adopted by Chelkak and Smirnov \cite{ChelkakSmirnov:ising} to prove conformal
invariance of the Ising model on bounded domains of the plane with Dobrushin
boundary conditions. The second approach was used by the authors
in \cite{BoutillierdeTiliere:iso_perio, BoutillierdeTiliere:iso_gen} to construct a Gibbs measure for the
critical Ising model on infinite isoradial graphs.

Before explaining these two approaches, let us mention some works by physicists on the $Z$-invariant model at and off criticality by Au-Yang and Perk \cite{Perk1,Perk2,Perk3} and by Martinez \cite{Martinez1,Martinez2}.

\subsection{Conformal invariance}

%

Let $\Omega$ be a bounded domain of $\mathbb{R}^2$, and fix two points $a$ and $b$ on the
boundary. For every $\varepsilon$, let $\GR_\varepsilon$ be an approximation
of $\Omega$ by a rhombus-graph\footnote{On the boundary, we put only half-rhombi
such that only ``black'' vertices are exposed on the boundary.} with rhombi of
side length $\varepsilon$. Let also  $a_\varepsilon$ and $b_\varepsilon$ be
approximations of $a$ and $b$, located at the center of half-rhombi on the
boundary.

Consider the Ising model on the isoradial graph $G_\varepsilon$ corresponding
to black vertices of $\GR_\varepsilon$, with Dobrushin boundary conditions, 
\emph{i.e} fixed spins $+1$ on one arc from $a_\varepsilon$ to $b_\varepsilon$, 
and $-1$ on the other. 
Introduce, following Chelkak and Smirnov \cite{ChelkakSmirnov:ising}, a twisted version of the partition
function, considered as a function on edges of $G_\eps$, or equivalently on inner rhombi of $\GR_\eps$. Let $z$
denote a generic rhombus of $\GR_\eps$,
\begin{equation*}
  \tilde{Z}_{G_\varepsilon,a_\varepsilon,z} = \tilde{Z}_\varepsilon(z)  = 
  \left(\sin \frac{\theta_z}{2}\right)^{-1} \sum_{\substack{\text{closed contours} \\ \text{ + curve } 
\gamma:a_\varepsilon \leadsto z}} 
  \left(\prod_{e: \text{ piece } \\ \text{of contour}} \tan\frac{\theta_e}{2}\right) e^{-\frac{i}{2} 
\text{wind}(\gamma)}
\end{equation*}
where $\text{wind}(\gamma)$ is the winding of the curve $\gamma$ from $a$ to $z$. 
If we remove the prefactor and the contribution of the winding, this would be
for $z=b_\varepsilon$ the low temperature expansion of the partition function
of the Ising model on $\Omega_\varepsilon$ with Dobrushin boundary conditions
between $a_\varepsilon$ and $b_\varepsilon$.

Now define for all inner rhombi $z$ of $\GR_\eps$,
\begin{equation*}
  F_{G_\varepsilon,a_\varepsilon,b_\varepsilon}(z)=F_{\varepsilon}(z) \propto \frac{\tilde{Z}_\varepsilon(z)}{\tilde{Z}_\varepsilon(b_\varepsilon)},
\end{equation*}
up to a multiplicative factor, introduced for technical reasons.

Chelkak and Smirnov prove that the function $F_{\varepsilon}$
and is discrete holomorphic by comparing configurations differing by local
arrangements near $z$, and that it solves some discrete Riemann-Hilbert
boundary value problem. 
They then deduce using approximation results \cite{ChelkakSmirnov:toolbox} that,
 as $\varepsilon\to 0$,
the function $F_{\varepsilon}$ converges to the function solving the analogue
continuous Riemann-Hilbert boundary value problem, which is conformally
covariant.

Moreover, they prove that this observable satisfies a martingale property with
respect to the growth of the curve from $a_\varepsilon$ to $b_\varepsilon$:
given the first $n$ steps of the curve $( \gamma_0 =a_\varepsilon,
\gamma_1,\dots, \gamma_n)$, the expected value of the observable
$F_{G_\varepsilon\setminus[\gamma_0,\gamma_{n+1}],\gamma_{n+1},b_\varepsilon}(z)$
is equal to
$F_{G_\varepsilon\setminus[\gamma_0,\gamma_{n}],\gamma_{n},b_\varepsilon}(z)$.

This martingale property together with the convergence to a conformally covariant object is 
sufficient to imply the convergence of the interface between $a$ and $b$ to a chordal SLE, 
with parameter $\kappa=3$.

In the same paper, Chelkak and Smirnov construct another discrete holomorphic observable using
not the loops separating clusters of spins, but those separating the clusters
in the corresponding random cluster model\footnote{also known as the
Fortuin-Kasteleyn percolation.} with $q=2$. This observable is a direct
generalization of the one introduced by Smirnov in \cite{Smirnov:isingZ2} for
the square lattice.
Again, this observable is discrete holomorphic and satisfies a martingale property. In this case, 
the scaling limit of the interface is a chordal SLE(16/3). 

\subsection{The two-dimensional Ising model as a dimer model}

\label{subseq:Ising-dimers}
It turns out that a general Ising model on a graph drawn on a surface without
boundary can be reformulated, through its contour expansion, as a dimer model on a decorated graph.
This correspondence is due to Fisher \cite{Fisher}. Since a lot of
exact computations can be carried out on the dimer model, this correspondence
has been proven to be a useful tool to study the Ising model
\cite{McCoyWu} (\cite{Kasteleyn}, \ldots)

The idea is to construct a \emph{decorated version} $\GD$ of the graph $G$ with
a measure preserving mapping between polygonal contours of $G$ and
dimer configurations of $\GD$. This is an example of \emph{holographic
reduction} \cite{Valiant}.  We now present the version used by the authors in
\cite{BoutillierdeTiliere:iso_perio,BoutillierdeTiliere:iso_gen}, which has the
advantage of using decoration with cyclic symmetry. For other versions of Fisher's correspondence, see
\cite{Fisher}, \cite{Cimasoni:KacWard}.

The decorated graph, on which the dimer configurations live, is constructed
from $G$ as follows. Every vertex of degree $k$ of $G$ is replaced by a {\em
decoration} consisting of $3k$ vertices: a triangle is attached to every edge
incident to this vertex, and these triangles are linked by edges in a circular
way, see Figure \ref{fig:decorated_graph} below. This new graph, denoted by
$\GD$, is also embedded on the surface without boundary and has vertices of
degree $3$. It is referred to as the {\em Fisher graph} of $G$. 

\begin{figure}[ht]
\begin{center}
\includegraphics[height=2cm]{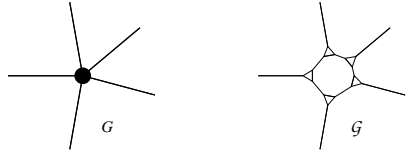}
\caption{Left: a vertex of $G$ with its incoming edges. Right: corresponding decoration in $\GD$.}
\label{fig:decorated_graph}
\end{center}
\end{figure}
%
Here comes the correspondence: to any contour configuration $\mathsf{C}$ coming
from the high-temperature expansion of the Ising model on $G$, we associate
$2^{|V(G)|}$ dimer configurations on $\GD$: edges present (resp. absent) in
$\mathsf{C}$ are absent (resp. present) in the corresponding dimer
configuration of $\GD$. Once the state of these edges is fixed, there is, for
every decorated vertex, exactly two ways to complete the configuration into a
dimer configuration. Figure \ref{fig:Fisher_correspondence} below gives an
example in the case where $G$ is the square lattice $\ZZ^2$.

\begin{figure}[ht]
\begin{center}
\includegraphics[width=\linewidth]{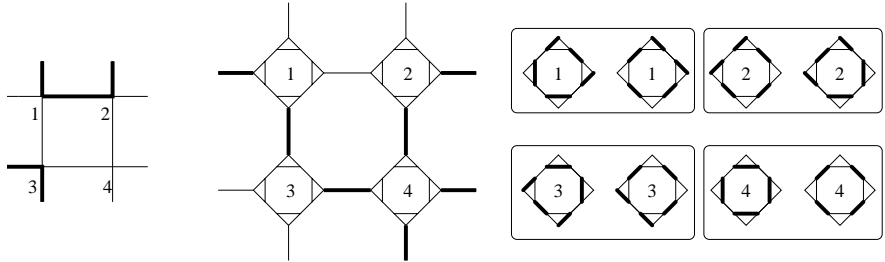}
\caption{Polygonal contour of $\ZZ^2$, and corresponding dimer configurations
  of the Fisher graph.}
\label{fig:Fisher_correspondence}
\end{center}
\end{figure}

In order to have a measure preserving mapping with dimers, the dimer weight function $\nu$ is given for 
the special case of the $Z$-invariant critical Ising model by,
\begin{equation*}
  \begin{cases}
    \nu_e = \cot\frac{\theta}{2} & \text{if $e$ comes from a rhombus of $G$ with half-angle $\theta$,} \\
    \nu_e = 1 & \text{if $e$ belongs to a decoration.}
  \end{cases}
\end{equation*}
%
%

A Kasteleyn matrix $K$ is constructed on $\GD$. One of the main result of \cite{BoutillierdeTiliere:iso_gen} is that
its inverse has a `local' expression in the spirit of that obtained by Kenyon in \cite{Kenyon3} for the inverse
Dirac operator, see also Theorem \ref{thm:dirac}. 
\begin{theorem}[\cite{BoutillierdeTiliere:iso_gen}]
If $x$ (resp. $y$) is a vertex of $\GD$ and belongs to the 
decoration corresponding to the vertex $\mathbf{x}$ (resp. $\mathbf{y}$) of $G$, then $(K^{-1})_{x,y}$ has the following 
expression:
\begin{equation*}
  (K^{-1})_{x,y} = \frac{1}{4\pi^2}\oint_{C_{x,y}} g_x(\lambda) g_y(-\lambda) 
\expo_{\mathbf{x}}(\mathbf{y}; \lambda) \log(\lambda)\mathrm{d}\lambda, 
\end{equation*}
where $g_x$ (resp. $g_y$) is a simple rational fraction of $\lambda$ that depends only on the geometry of 
the decorated graph in an immediate neighborhood of $x$ (resp. of $y$), and the contour of integration 
$C_{x,y}$ is a simple closed curve containing all poles of the
integrand, and avoiding the half-line starting from $x$ in the direction of $y$.
\end{theorem}

Then, as for the dimer model, we obtain an explicit expression for a Gibbs measure, and recover the explicit formula
for the partition function, obtained by Baxter \cite{Baxter:Zinv}, at the critical point, see \cite{BoutillierdeTiliere:iso_gen} for precise statements. 

\section{Other models}

We briefly discuss now some aspects of other models of statistical mechanics: first some models related to the Laplacian, then the $q$-Potts model and its random cluster representation, and finally the $6$-vertex and $8$-vertex models.

\subsection{Random walk and the Green function}

We already mentioned the Laplace operator on functions of vertices, corresponding to conductances on edges given by $c(e)=c(\theta_e)=\tan\theta_e$.
For this particular choice of conductances, the associated random walk is a martingale and the covariance matrix 
associated to a step is scalar. Up to a time reparametrization, its scaling limit is a standard two-dimensional 
Brownian motion \cite{ChelkakSmirnov:toolbox, Beffara}.
Kenyon \cite{Kenyon3} proves a \emph{local formula} for the inverse of the Laplacian, the  Green function $\Delta^{-1}$, on an infinite isoradial graph:
\begin{equation*}
  \forall x,y \in V(G), \quad \Delta^{-1}_{x,y} = - \frac{1}{8i\pi^2} \oint_C \expo_{x}(y; \lambda) \frac{\log \lambda}{\lambda} \mathrm{d}\lambda.
\end{equation*}
where  $C$ is a positively oriented contour containing in its interior all the poles of the discrete exponential 
function, and the cut of the determination of the logarithm. He also proves precise asymptotics of $\mathrm{G}$, improved in \cite{Buecking}.
This Green function gives information on the random walk on $G$ with conductances $\tan\theta$,
 but can also be used to gather properties of other models from statistical mechanics with an interpretation in terms of electric network, such as random spanning trees.

A \emph{spanning tree} of $G$ is an acyclic connected subgraph of $G$ whose vertex set contains all vertices of $G$. 
On a finite graph, if the weight of a spanning tree equals the product of the conductances of its edges, 
the partition function is given, via Kirchhoff's matrix-tree theorem \cite{Kirchhoff:1847}, by the determinant of 
any principal minor of the Laplacian. One can construct a measure on spanning trees of $G$ as limit of measures on 
finite graphs where a spanning tree would have a probability proportional to its weight. The statistics of edges 
present in the random spanning tree are given by a determinantal process on edges: the probability
that edges $\{e_1=v_1w_1$, \ldots, $e_k=v_kw_k\}$ are present is given by
\begin{equation*}
  P_{\text{tree}}(\{e_1,..,e_k\}) = \left(\prod_{i=1}^k \tan \theta_i\right) \det_{1\leq i,j\leq k}\bigl( H(e_i,e_j) \bigr),
\end{equation*}
where  $H$ is the impedance transfer matrix \cite{BurtonPemantle} defined by
\begin{equation*}
  H(e_i,e_j)=\Delta^{-1}(v_i,v_j)-\Delta^{-1}(v_i,w_j)-\Delta^{-1}(w_i,v_j)+\Delta^{-1}(w_i,w_j).
\end{equation*}

\subsection{$q$-Potts models and the random cluster model}

A natural generalization of the Ising model is the Potts model with $q$ states (or colors). Spins at neighboring vertices interact only if they have different colors. The energy of a spin configuration $\sigma$ is $
  \mathcal{E}(\sigma)=-\sum_{e=uv\in E} J_e \mathbb{I}_{\sigma_u\neq\sigma_v}$.
  For $q=2$, we recover the Ising model. This model can be mapped to the 
(inhomogeneous) random cluster model \cite{Grimmett:FK} as follows: 
\begin{itemize}
  \item spins of different colors are in different percolation clusters,
  \item neighboring spins sharing an edge $e$ of the same color are connected with probability, $p_e=1-e^{-2J_e}$.
\end{itemize}
Once again, we want the interaction constants to be functions of the half-angle of the corresponding rhombus. It is possible to take them to satisfy the Yang-Baxter equation. Moreover, if we impose a generalized self-duality, then there is a unique choice for the function $p$ \cite{Kenyon:trieste}:
\begin{equation*}
  p(\theta)=\frac{\sqrt{q} \sin\bigl(\frac{2r}{\pi}\theta\bigr)}{\sin\bigl(\frac{2r}{\pi}(\frac{\pi}{2}-\theta)\bigr) +\sqrt{q}\sin\bigl(\frac{2r}{\pi}\theta\bigr)}, \quad \text{with }r=\cos^{-1}\left(\sqrt{q}/2 \right)
\end{equation*}
Although a formula for the free energy of that model is known \cite{Baxter:exactly}, little is known except for the case $q=2$ corresponding to the Ising model.

\subsection{6-vertex and 8-vertex models}
To conclude, let us briefly mention that the 6- and 8-vertex models, introduced originally on the square lattice, 
have natural generalization on isoradial graphs, or more precisely on the dual of its rhombic graph ${\GR}^*$, which is a $4$-regular graph.
A configuration of the 6-vertex model (resp. of the 8-vertex model) is an orientation of the edges of ${\GR}^*$ such 
that the number of ingoing edge (and thus of outgoing edges) at each vertex is equal to $2$ (resp. is even).
These models can be solved through Bethe Ansatz on the square lattice (in the sense that their free energy can be 
computed \cite{Lieb,Baxter:8V}), and possess a very rich algebraic structure. There are some studies on their 
$Z$-invariant generalizations \cite{Baxter:perimeter,Baxter:Zinv}, but many questions still need to solved.

\bibliographystyle{alpha}
\bibliography{survey}
\end{document}